\documentclass[10pt, reqno]{amsart}

\usepackage{latexsym}
\usepackage{amsthm}
\usepackage{amsmath}
\usepackage{amsfonts}
\usepackage{euscript}
\usepackage[all, cmtip]{xy}
\usepackage{latexsym,tabularx}
\usepackage{amscd}
\usepackage{graphics}
\usepackage{amssymb}
\usepackage{epsfig}
\usepackage{fancyhdr}

\newtheorem{lem}{Lemma}[section]

\newtheorem{defn}{Definition}[section]

\newtheorem{prop}{Proposition}[section]

\begin{document}
\title[Rank two vector bundles with canonical determinant]{Remarks on Ruled surfaces and rank two bundles with canonical determinant and 4 sections.}
\author{Abel Castorena and Graciela Reyes-Ahumada.}
\address{Centro de Ciencias Matem\'aticas. Universidad Nacional Aut\'onoma de M\'exico, Campus Morelia\\
Apdo. Postal 61-3(Xangari). C.P. 58089\\
Morelia, Michoac\'an. M\'exico.}
\email{abel@matmor.unam.mx, grace@matmor.unam.mx}
\begin{abstract} Let $C$ be a smooth irreducible complex projective curve of genus $g$ and let $B^k(2,K_C)$ be the Brill-Noether locus parametrizing classes of (semi)-stable vector bundles $E$ of rank two with canonical determinant over $C$ with $h^0(C,E)\geq k$. We show that $B^4(2,K_C)$ has an irreducible component $\mathcal B$ of dimension $3g-13$ on a general curve $C$ of genus $g\geq 8$. Moreover, we show that for the general element $[E]$ of $\mathcal B$, $E$ fits into an exact sequence $0\to\mathcal O_C(D)\to E\to K_C(-D)\to 0$ with $D$ a general effective divisor of degree three, and the corresponding coboundary map $\partial: H^0(C,K_C(-D))\to H^1(C,\mathcal O_C(D))$ has cokernel of dimension three.
\end{abstract}
\maketitle
\section{Introduction and statement of the result.}

Let $C$ be a smooth irreducible complex projective curve of genus $g$ and let $U_C(r,d)$ be the moduli space of (semi)stable vector bundles of rank $r$ and degree $d$ on $C$. Inside $U_C(r,d)$, consider the Brill-Noether locus $B^k(r,d)$ parametrizing classes of vector bundles $[E]\in U_C(r,d)$ having at least $k$ linearly independent sections. Because of the interest geometry behind such loci, higher rank Brill-Noether theory on algebraic curves has been extensively studied and it is actually an active research area in Algebraic Geometry; see e.g. (\cite{IT}), for an overview of some main results, and references below, for more recent ones. Despite these facts, several basic questions concerning e.g. non-emptiness, dimensionality, irreducibility, etc., are still open. Let $B(2,K_C)\subset U_C(2, 2g-2)$ be the scheme which parametrizes classes of (semi)stable rank-two vector bundles $E$ with $\text{det}(E)=\wedge^2E=K_C$; this scheme is defined as the fiber at $K_C$ of the determinant map $\wedge^2: U_C(2,2g-2)\to\text{Pic}^{2g-2}( C), E\to\wedge^2 E$. $B(2,K_C)$ is smooth and irreducible of dimension $3g-3$. Inside $B(2,K_C)$, for any integer $k\geq 0$ define the Brill-Noether locus $B^k(2,K_C):=\{[E]\in B(2,K_C):h^0(C,E)\geq k\}$. These loci have a scheme structure as suitable degeneracy loci and it is known that the expected dimension of $B^k(2,K_C)$ is given by the Brill-Noether number $\rho_{K_C}(2,k,g):=3g-3-\binom{k+1}{2}$ (see \cite{BF}, \cite{M1}). For a description of $B^k(2,K_C)$ for low genus we refer to e.g. (\cite{BF}, \cite{M2}, \cite{M3}). For $[E]\in B^k(2,K_C)$, the infinitesimal behavior of $B^k(2,K_C)$ at the point $[E]$ is governed by the symmetric Petri map $P_E:\text{Sym}^2(H^0(C,E))\to H^0(C,\text{Sym}^2(E))$, that is, the tangent space to $B^k(2,K_C)$ at $[E]$ is identified with the orthogonal to the image of $P_E$. In (\cite{T2}) it has been proved the injectivity of the symmetric Petri map $P_E$ on a general curve of genus $g\geq 1$, when $B^k(2,K_C)$ is assumed to be non-empty. A different approach for the injectivity of the symmetric Petri map for $g\leq 9$ and $k<7$ is given in (\cite{B}). The injectivity of the symmetric Petri map implies that on a general curve $C$, components of the right dimension in $B^k(2,K_C)$ are smooth. 

\vskip1mm

\noindent A general result on non-emptiness and existence of components in $B^k(2,K_C)$ of the right dimension on a general curve of genus sufficiently large is given in ({\cite{T1}}), where the proof uses the theory of limit linear series for higher rank. For generalizations of this result and others results about non-emptiness, irreducibility and smoothness, we refer to e.g (\cite{CF2}, \cite{LNP}, \cite{O}, \cite{VGI}, \cite{Z}). Focusing on (\cite{CF2}), as a by-product of the more general approach developed therein, the authors prove in particular the existence of irreducible components of $B^k(2,K_C)$ of the right dimension, for $C$ a general curve of genus $g >3$ and for $k \leq 3$; this has been done by studying certain determinantal loci inside the space of extensions 
$\text{Ext}^1(K_C(-D_{k-1}),\mathcal O_C(D_{k-1}))$, for a general effective divisor $D_{k-1}$ of degree $k-1\leq 2$. For a general curve of genus $g\geq 5$, the existence of irreducible components in $B^4(2, K_C)$ of the right dimension is a particular case of Theorem 1.1 given in (\cite{T1}). We follow the approach as in (\cite{CF2}) to show that there exists an irreducible component $\mathcal B\subseteq B^4(2,K_C)$ of (expected) dimension $\rho_{K_C}(2,4,g)$ on a general curve of genus $g\geq 8$ moreover we show that the general vector bundle in $\mathcal B$ is given as an extension of suitable line bundles as in (\cite{CF2}), to which the reader is referred for more details. In order to state our main theorem we recall here some definitions and notation.  

\vskip2mm

\noindent Let $C$ be a non-hyperelliptic curve of genus $g\geq 3$ and let $d$ be an integer such that $2g-2\leq d\leq 4g-4$. Let 
$\delta\leq d$ be a positive integer. Let $L\in\text{Pic}^{\delta}( C), N\in\text{Pic}^{d-\delta}( C)$ be line bundles on $C$; the space 
$\text{Ext}^1(L,N)$ parametrizes isomorphism classes of extensions 

$$0\to N\to E \to L\to 0.$$ 

Any $u\in\text{Ext}^1(L,N)$ gives rise to a degree $d$ rank-two vector bundle $E_u$ that fits in an exact sequence 

$$(u):\hskip5mm 0\to N\to E_u \to L\to 0.$$ 

\noindent To get a (semi)stable vector bundle, a necessary condition is $2\delta-d\geq 0$ (cf. \cite{CF2}, Remark 2.5 and Section 4.1). By Riemann-Roch, $\text{Ext}^1(L,N)\simeq H^1(C,N\otimes L^{\vee})\simeq H^0(C,K_C\otimes L\otimes N^{\vee})^{\vee}$, then $m:=\text{dim}(\text{Ext}^1(L,N))=2\delta-d+g-1$ if $L\ncong N$, and $m=g$ when $L\simeq N$. 

\vskip1mm

\noindent Suppose that $N$ is a special line bundle. For any $u\in\text{Ext}^1(L,N)$, consider the coboundary map induced by cohomology in the exact sequence $(u)$: 

$$\partial_u: H^0(L)\to H^1(N).$$ 

\noindent By exactness in $(u)$ we have that $h^1(C,E_u)=h^1(L)+\text{dim}(\text{Coker}(\partial_u))$. For any integer $t>0$, consider the degeneracy locus 

$$\mathcal W_t:=\{u\in\text{Ext}^1(L,N):\text{dim(Coker)}(\partial_u)\geq t.\}, $$ 

\noindent which has a natural description of determinantal scheme. If $m>0$ and $\mathcal W_t\neq\emptyset$, any irreducible component $\Delta_t\subseteq\mathcal W_t$ is such that $\text{dim}(\Delta_t )\geq\text{min}\{m,m-t(t+h^0(L)-h^1(N))\}$, where the right-hand-side is the expected dimension (cf. \cite{CF2}, Section 5.2). 

\noindent For the particular case of $B^4(2,K_C)$, if one tries to construct stable vector bundles with canonical determinant by extensions with $h^1(C,N)=0$ (so we are forced to have $h^0(C,L)=4$), we can get non-emptiness of $B^4(2,K_C)$ but we cannot describe the general element of a component, so the condition on 
$N$ to be a special line bundle give us an easier way to study this problem. Our main result is the following:

\vskip2mm

\noindent{\bf{Theorem.}} Let $C$ be a curve of genus $g\geq 8$ with general moduli. Then $B^4(2,K_C)\neq\emptyset$. Moreover, there exists an irreducible component 
$\mathcal B\subseteq B^4(2,K_C)$ of the expected dimension $\rho_{K_C}(2,4,g)=3g-13$ and whose general point $[E]$ fits into an exact sequence

$$0\to \mathcal O_C(D)\to E\to K_C(-D)\to 0,$$

\vskip1mm

where,

\vskip1mm

\noindent (i).- $D$ is a general effective divisor of degree 3, and

\vskip1mm

\noindent (ii).- $E=E_u$ with $u\in\Delta_3\subseteq\mathcal W_3\subset\text{Ext}^1(K_C(-D),\mathcal O_C(D))$ general in $\Delta_3$, where $\Delta_3$ is an irreducible component of dimension $3g-15$, whose general element $u\in\Delta_3$ satisfies 
$\text{dim(Cokernel}(\partial_u))=3$.

\vskip1mm

\noindent{\bf{Acknowledgements.}} The authors give thanks to Ciro Ciliberto, Flaminio Flamini, Elham Izadi and Christian Pauly for fruitful discussions and suggestions.

\noindent The first named author was partially supported by sabbatical fellowships from CONACyT (project 133228) and from PASPA (DGAPA, UNAM-M\'exico). The second named author was supported by a CONACyT scholarship and partially supported by CONACyT research project 166158 and LAISLA-M\'exico. 

\section{Preliminaries.}
We remind some definitions and results which are used for the proof of our main theorem. For more details on the topics of this section we refer the reader to e.g. (\cite{ACGH}, \cite{CF2}, \cite{H}, \cite{S}). 

\vskip1mm

\subsection{Ruled surfaces and sections.} Let $E$ be a vector bundle over a smooth irreducible complex projective curve $C$ of genus $g$. The speciality of $E$ is defined as $i(E)=h^1(C,E)$. $E$ is said special if $i(E)>0$. 
\vskip1mm

\noindent If the rank of $E$ is two, let $S:=\Bbb P(E)$ be the (geometrically) ruled surface associated to it, with structure map $p:S\to C$.  For any $x\in C$ we denote by $f_x=p^{-1}(x)\simeq\Bbb P^1$. We denote by $f$ a general fiber of $p$ and by $\mathcal O_S(1)$ the tautological line bundle on $S$. We write $d=\text{deg}(E)=\deg(\wedge^2 E)=\text{deg (det}(E))$. There is a map $s:C\to S$ such that $p\circ s=Id_C$ whose image we denote by $H$ and such that $\mathcal O_S(H)=\mathcal O_S(1)$. Assume that $h^0(C,E)>0$; thus an element in the linear system $|\mathcal O_S(1)|$ is denoted by $H$. For any $D\in\text{Div}(C)$ ( or in 
$\text{Pic}( C)$ ) we put $f_D:=p^*(D)$. If $\Gamma$ is a divisor on $S$, we set $\text{deg}(\Gamma):=H\cdot\Gamma$. We have that $d=\text{deg}(E)=H^2=\text{deg}(H)$. We recall that $\text{Pic}(S)\simeq\Bbb Z[\mathcal O_S(1)]\oplus p^*(\text{Pic} (C))$, moreover $\text{Num}(S)\simeq\Bbb Z\oplus\Bbb Z$ and is generated by the classes $H$, $f$ satisfying $H\cdot f=1, f^2=0$ (cf. \cite{H}, chapter V). Any element of $\text{Pic}(S)$ corresponds to a divisor on $S$ of the form $nH+f_B, n\in\Bbb Z, B\in\text{Div}(C)$, as an element of $\text{Num}(F)$, corresponds to $nH+bf, b=\text{deg}(B)$. For any $n\geq 0$ and for any $D\in\text{Div}(C)$, the linear system $|nH+f_D|$ if non-empty, is said to be {\it{$n$-secant }} to the fibration $p:S\to C$ since its general elements meets $f$ at $n$ points. An element $\Gamma\in |H+f_D|$ is called {\it{unisecant curve}} of $S$ (or to the fibration $p:S\to C$). The irreducible unisecants of $S$ are smooth and isomorphic to $C$ and are called {\it{sections of $S$}}. We denote by $\sim$ linear equivalence of divisors and by $\equiv$ numerical equivalence of divisors. 

\vskip1mm

\noindent We recall that there is a one-to-one correspondence between sections $\Gamma$ of $S$ and surjective maps $E\to\!\!\to L$ with $L$ a line bundle on $C$ (cf. \cite{H}, chapter V), then one has an exact sequence $0\to N\to E\to L\to 0$, where $N\in\text{Pic}( C)$. The surjection $E\to\!\!\to L$ induces an inclusion $\Gamma=\Bbb P(L)\subset S=\Bbb P(E)$. If $L=\mathcal O_C(B), B\in\text{Div}( C)$ with $b=\text{deg}(B)$, then $b=H\cdot\Gamma$ and $\Gamma\sim H+f_N$ where $N=L\otimes\text{det}(E)^{\vee}\in\text{Pic} (C)$. For a section $\Gamma$ that corresponds to the exact sequence $0\to N\to E\to L\to 0$, we have that $\mathcal N_{\Gamma/S}$, the normal sheaf of $\Gamma$ in $S$, is such that $\mathcal N_{\Gamma/S}=\mathcal O_{\Gamma}(\Gamma)\simeq N^{\vee}\otimes L$. In particular $\text{deg}(\Gamma)=\Gamma^2=\text{deg}(L)-\text{deg}(N)$. If $\Gamma\sim H-f_D$, i.e. when $N=\mathcal O_C(D)$, then $|\mathcal O_S(\Gamma)|\simeq \Bbb P(H^0(C,E(-D)))$. For $\Gamma\in\text{Div}(S)$, we define $\mathcal O_{\Gamma}(1):=\mathcal O_S(1)\otimes\mathcal O_{\Gamma}$. The {\it{speciality of $\Gamma$}} is defined as $i(\Gamma):=h^1(\Gamma,\mathcal O_{\Gamma}(1))$. {\it{$\Gamma$ is said special if $i(\Gamma)>0$}}. 

\vskip1mm

\subsection{Hilbert scheme of unisecant curves and Quot-Scheme.} For any $n\in\Bbb N$, denote by $\text{Div}^{1,n}(S)$ the Hilbert scheme of unisecant curves of $S$ which are of degree $n$ with respect to $\mathcal O_S(1)$. Since elements of $\text{Div}^{1,n}(S)$ correspond to quotients of $E$, therefore $\text{Div}^{1,n}(S)$ can be endowed with a natural structure of Quot-Scheme (cf. \cite{S}, section 4.4), and one has an isomorphism (see e.g. \cite{CF2}, Section 2.4.) 

$$\Phi_{1,n}:\text{Div}^{1,n}(S)\xrightarrow{\simeq}\text{Quot}^C_{E,n+a-g+1},\hskip5mm \Gamma\to\{E\to\!\!\to L\oplus\mathcal O_A\},$$

\noindent where $\Gamma$ is the unisecant $\Gamma'+f_A$ with $A\in\text{Div}(C)$, and $\Gamma'$ is the section corresponding to $E\to\!\!\to L$. The morphism $\Phi_{1,n}$ allows to identify tangent spaces $H^0(\Gamma,\mathcal N_{\Gamma/S})\simeq T_{[\Gamma]}(\text{Div}^{1,n}(S))\simeq\text{Hom}(N\otimes\mathcal O_C(-A), L\oplus\mathcal O_A)$ and obstruction spaces $H^1(\Gamma,\mathcal N_{\Gamma/S})\simeq\text{Ext}^1(N\otimes \mathcal O_C(-A),L\oplus\mathcal O_A)$ (cf. \cite{S}, Section 4.4). 

 \vskip1mm

\noindent For any integer $m$ such that $m\leq\overline m:=\lfloor\frac{d-g+1}{2}\rfloor$, one can consider also the set $Q_m(E):=\{N\subset E: N\text{ is an invertible subsheaf of }E,\text{deg}(N)=m\}$. This set has a natural structure of Quot-scheme.  Let $\Gamma$ be any unisecant curve on $S$ of degree $n$, corresponding to the exact sequence $0\to N\to E\to L\oplus\mathcal O_A\to 0$, where $L, N$ are line bundles such that $n=\text{deg}(L)+\text{deg}(A)$ and $m:=\text{deg}(N)=d-n$. One has an isomorphism 
$$\text{Div}^{1,n}(S)\xrightarrow{\Phi_{n,m}}Q_m(E),\hskip1mm \Phi_{n,m}([\Gamma])=[N].$$ 
\noindent Moreover, one has a natural morphism $\pi_m: Q_m(E)\to\text{Pic}^m(C), N\to[N]$, and we denote $W_m(E):=\text{Im}(\pi_m)$. For any $[N]\in W_m(E)$, one has that $\pi_m^{-1}([N])\simeq \Bbb P(H^0(E\otimes N^{\vee}))$ (see \cite{Gh}, p. 199), this implies in particular that $\pi_m$ has connected fibres. For every integer $p\geq 0$, one can consider the loci $Q_m^p(E):=\{N\subset Q_m(E): h^0(E\otimes N^{\vee})\geq p+1\}$, $W_m^p(E):=\{N\subset W_m(E): h^0(E\otimes N^{\vee})\geq p+1\}$. For some properties of these loci we refer the reader to e.g. (\cite{Gh}).

\vskip4mm
 
By the Quot-scheme structure of $\text{Div}^{1,n}(S)$, the universal quotient $\mathcal Q_{1,\delta}\xrightarrow{\pi}\text{Div}^{1,n}(S)$ gives a morphism 
$\text{Proj}(\mathcal Q_{1,\delta})\xrightarrow{\pi}\text{Div}^{1,n}(S)$. We define 
 $$\mathcal S^{1,n}:=\{\Gamma\in\text{Div}^{1,n}(S): R^1\pi_*(\mathcal O_{\mathcal Q_{1,\delta}}(1))_{\Gamma}\neq 0\}.$$ 
 
\noindent This locus supports the {\it{scheme that parametrizes degree $n$, special unisecant curves of $S$}}.

\vskip1mm

\begin{defn}( cf. \cite{CF2}, Def. 2.2, 2.10) Let $\Gamma\in\text{Div}^{1,n}(S)$ be. We say that:

(i).- $\Gamma$ is linearly isolated if $\text{dim}(|\mathcal O_S(\Gamma)|)=0$;

(ii).- $\Gamma$ is algebraically isolated if $\text{dim}(\text{Div}^{1,n}(S))=0$.

\vskip1mm

\noindent Let $\Gamma$ be a special unisecant of $S$. Assume that $\Gamma\in\mathcal F$, where $\mathcal F\subset\text{Div}^{1,n}(S)$ is a subscheme.

(iii).- $\Gamma$ is specially unique in $\mathcal F$, if $\Gamma$ is the only special unisecant in $\mathcal F$.

(iv)-  $\Gamma$ is specially isolated in $\mathcal F$, if $\text{dim}_{\Gamma}(\mathcal F\cap\mathcal S^{1,n})=0$.
\end{defn}

\noindent When $\mathcal F=|\mathcal O_S(\Gamma)|$, $\Gamma$ is said to be {\it{linearly specially unique}} in case (iii) and {\it{linearly specially isolated}} in case (iv). 

\noindent When $\mathcal F=\text{Div}^{1,n}(S)$, $\Gamma$ is said to be {\it{algebraically specially unique}} in case (iii) and {\it{algebraically specially isolated}} in case (iv). 

\subsection{A result on deformation theory. }Let $Y$ be a smooth projective variety and $\jmath: X\subset Y$ be a closed smooth subvariety. Let $\mathcal I_X\subset\mathcal O_Y$ be the ideal sheaf of $X$. We have an inclusion of tangent sheaves $T_X\subset T_Y|_X$ and a restriction morphism $T_Y\to\!\!\to T_Y|_X$. Let $T_Y<X>\subset T_Y$ be the image inverse of $T_X$ under the restriction morphism. The sheaf $T_Y<X>$ is called the {\it{sheaf of germs of tangent vectors to $Y$ which are tangents to $X$}}. This is a coherent sheaf of rank $\text{dim} (Y)$ on $Y$. We have a restriction map $R:T_Y<X>\to\!\!\to T_X$ giving the exact sequence $0\to T_Y(-X)\to T_Y<X>\to T_X\to 0,$ where $T_Y(-X)$ is the vector bundle of tangent vectors of $Y$ vanishing along $X$. Let $H^1(Y,T_Y<X>)\xrightarrow{H^1(R)}H^1(X, T_X)$ be the map in cohomology induced by the above exact sequence. The map $H^1( R)$  associates to a first-order deformation of $(Y,X)$ the corresponding first-order deformation of $X$. We have (cf. \cite{S}, Proposition 3.4.17):

\vskip1mm

\noindent{\it{First-order deformations of the pair $(Y,X)$. }}The infinitesimal deformations of the pair $(Y,X)$ (equivalently of the closed embedding $\jmath$) are controlled by the sheaf 
$T_Y<X>$, that is,  

\vskip1mm

\noindent $(i)$.- The obstructions lie in $H^2(Y,T_Y<X>)$.

\noindent $(ii)$.- First-order deformations are parametrized by $H^1(Y,T_Y<X>)$ and the space $H^0(Y,T_Y<X>)$ parametrizes infinitesimal automorphisms.

\subsection{Extensions of line bundles and the Segre invariant.} Here we follow (\cite{CF2} Sections 2.1 and 4.1), to which we refer the reader for more details. 
\vskip1mm

\noindent For a rank two vector bundle $E$, the {\it{Segre invariant}} $s(E)$ of $E$ is defined as $s(E):=\text{deg}(E)-2(\text{max}\{\text{deg}(N)\})$, where the maximum is taken among all sub-line bundles $N$ of $E$. The bundle $E$ is stable (resp. semi-stable) if $s(E)>0$ (resp $s(E)\geq 0$). For any $A\in\text{Pic}( C)$, one has $s(E)=s(E\otimes A)$. 

\vskip1mm

\noindent For $\delta\leq d$ positive integers, take $L\in\text{Pic}^{\delta}( C)$ special and effective, whereas $N\in\text{Pic}^{d-\delta}( C)$ arbitrary. As in the Introduction, any $u\in\text{Ext}^1(L,N)$ gives rise to an exact sequence $(u): 0\to N\to E=E_u\to L\to 0$, which therefore corresponds to a point in $\Bbb P:=\Bbb P(\text{Ext}^1(L,N))$. Note that the bundle $\mathcal E_e:=E\otimes N^{\vee}$ fits into an exact sequence $(e): 0\to \mathcal O_C\to \mathcal E_e\to K_C\otimes N^{\vee}\to 0$. 

\vskip1mm

\noindent When $\text{deg}(L\otimes N^{\vee})\geq 2$,  $\text{dim}(\Bbb P)\geq 3$, then the map $\phi: C\to\Bbb P$ induced by $|K_C\otimes L\otimes N^{\vee}|$ is a morphism and posing $X:=\phi(C)\subset\Bbb P$, for any positive integer $h$ one can consider 
$\text{Sec}^h(X)$,  the $h^{th}$-secant variety of $X$, defined as the closure of the union of all linear spaces $<\phi(D)>\subset\Bbb P$ for general $D\in\text{Sym}^{(h)}(C)$. One has the following

\vskip1mm

\begin{prop} Let $2\delta-d\geq 2$. For any integer $\sigma$ such that $\sigma \equiv 2\delta-d(mod\ 2)$ and 
$4+d-2\delta\leq \sigma\leq 2\delta-d$, one has 
$$s(\mathcal E_e)\geq \sigma \Leftrightarrow e\notin Sec_{\frac{1}{2}(2\delta-d+\sigma-2)}(X).$$
\end{prop}

\noindent See (\cite{LN}), Proposition 1.1.

\vskip1mm

\noindent The previous result gives in particular geometric conditions for $E_u$ as above to be 
(semi)stable. When $N$ is assumed to be special, as in \cite{CF2}, Section 5.2, one can define for any positive integer 
$t\geq 1$ the locus $\mathcal W_t=\{u\in\text{Ext}^1(L,N) | \partial_u:H^0(L)\to H^1(N)\text{ has cokernel of dimension }\geq t \}$, which has a natural structure of determinantal scheme, as such of expected dimension $m-t(t+h^0(L)-h^1(N))$, where $m=\text{dim}(\text{Ext}^1(L,N))$ as in Introduction. 

\begin{defn} (cf. \cite{CF2}, Def. 5.12). For $h^0(L)\geq h^1(N)\geq t\geq 1$ integers, assume that:

\vskip1mm

(i).- There exists an irreducible component $\Delta_t\subset\mathcal W_t$ of the expected dimension $\text{dim}(\Delta_t)=m-t(t+h^0(L)-h^1(N));$

\vskip1mm

(ii).- For $u\in\Delta_t$ general, $\text{dim(Coker)}(\partial _u)=t$.

\noindent Any such $\Delta_t$ is called a {\it{Good component}} of $\mathcal W_t$. See (\cite{CF2}, definition 5.12).\end{defn}

\noindent Theorems 5.8, 5.17 in (\cite{CF2}) give sufficient conditions for the existence of good components. These have been given by the use of suitable multiplication maps and an alternative description of 
$\mathcal W_t$. Precisely, for any subspace $W\subset H^0(K_C\otimes N^{\vee})$ of dimension at least $t$, which turns out to be 
$W=(\text{Coker}(\partial_u))^{\vee}$, one considers the natural multiplication map of sections $\mu:H^0(K_C\otimes N^{\vee})\otimes H^0(L)\to H^0(K_C\otimes L\otimes N^{\vee})$. Let 

$$\mu_W: W\otimes H^0(L)\to H^0(K_C\otimes L\otimes N^{\vee})\hskip22mm(+)$$ 

\noindent be the restriction of $\mu$ to $W$; then $\mathcal W_t$ is (cf. \cite{CF2}, Remark 5.7 for details): 

\vskip2mm

\noindent $\mathcal W_t=\{u\in H^0(K_C\otimes L\otimes N^{\vee})^{\vee}:\exists W\subset H^0(K_C\otimes N^{\vee})\text{ such that dim}(W)\geq t\text{ and Im}(\mu_W)\subset\{u=0\}\}.\hskip62mm (\star)$ 

\vskip3mm

\noindent This is the description of $\mathcal W_t$ we will use later on (cf. Lemma 3.2 and Step 1 of the proof of the main theorem).

\vskip1mm

\subsection{Global space of extensions.}

\noindent For the general case of the following construction we refer to e.g (\cite{ACGH}, \cite{CF2} Section 6.). 

\vskip1mm

\noindent Let $C$ be a general curve of genus $g\geq 3$. Let $0\to N\to E\to L\to 0$ be an exact sequence such that $\text{deg}(N)=d-\delta>0, \text{deg}(L)=\delta>0$ and $h^0(L)\cdot h^1(L)>0,h^0(N)\cdot h^1(N)>0$. Let $\ell:=h^0(L), j:= h^1(L), n:=h^0(N), r:=h^1(N)$. Set $\mathcal Y:=\text{Pic}^{d-\delta}( C)\times W^{\ell-1}_{\delta}( C),\hskip2mm Z:=W^{n-1}_{d-\delta}( C)\times W^{\ell-1}_{\delta} (C)\subset\mathcal Y$. One has $\text{dim}(\mathcal Y)=g+\rho_L$ and 
$\text{dim}(Z)=\rho_L+\rho_N$, where for any line bundle $M$ we denote by $\rho_M$ (or simply $\rho$, when $M$ is understood) the usual Brill-Noether number.

\noindent Since $C$ is general, we moreover have that $\mathcal Y$ and $Z$ are irreducible when $\rho>0$; otherwise one replace the (reducible) zero-dimensional Brill-Noether locus with one of its irreducible component to construct $Y$ and $Z$ as above.
\vskip1mm

\noindent When $2\delta-d\geq 1$ one can construct a vector bundle $\mathcal E\to\mathcal  Y$ of rank $m=2\delta-d+g-1$ (see \cite{ACGH}, p. 176-180, \cite{CF2} Section 6.), together with a projective bundle morphism $\gamma: \Bbb P(\mathcal E)\to\mathcal Y$ where, for $y=(N,L)\in\mathcal Y$, the fiber $\gamma^{-1}(y)=\Bbb P(\text{Ext}^1(N,L))=\Bbb P$. We have that $\text{dim}(\Bbb P(\mathcal E))=\text{dim}(\mathcal Y)+m-1,\text{ and }\text{dim}(\Bbb P(\mathcal E)|_Z)=\text{dim}(Z)+m-1$. Since (semi)stability is an open condition, for $2\delta-d\geq 2$ there is an open, dense subset $\Bbb P(\mathcal E)^0\subseteq\Bbb P(\mathcal E)$ and a morphism $\pi_{d,\delta}:\Bbb P(\mathcal E)^0\to U_C(d)$. 

\vskip1mm

\noindent In (\cite{CF2}, Sections 6,7) the authors study the image and fibers of the map $\pi_{d,\delta}$ under certain numerical conditions on $d,\delta$ and in a more general context with respect to the line bundles $L$ and $N$. They give in a different way a proof on the existence of irreducible and regular components in Brill-Noether loci $B^j(r,d)$ and $B^k(2,K_C), k\leq 3$. Some of such components are the (dominant) image under $\pi_{d,\delta}$ of certain degeneracy loci in $\Bbb P(\mathcal E)$ satisfying similar conditions to the good components of Definition 2.2. Such loci are called by the authors {\it{Total good components}} (see \cite {CF2}, Definition 6.13). Following this approach on a general curve $C$ of genus $g\geq 8$, for $k=4$ we obtain an irreducible (not good) component $\Bbb P(\widetilde{\Delta_3})\subset\Bbb P(\mathcal E)$ of dimension $3g-13$ that satisfies the conditions of the main theorem and fill-up an irreducible component 
$\mathcal B\subseteq B^4(2,K_C)$ of dimension $3g-13$. We prove this in the next section. 

\section{Proof of the Theorem.} 

For the following Lemma we adapt an argument of Lazarsfeld (cf. \cite{G}, Theorem 1.1). 

\vskip2mm

\begin{lem} Let $C$ be a non-hyperelliptic curve of genus $g\geq 8$ and let $D=q_1+q_2+q_3$ be a general effective divisor of degree 3 on $C$. There exists a rank two vector bundle $F$ on $C$ with the following properties: 

\vskip1mm

(i).- $\text{det}(F)=K_C(-D)$, $h^0(F)=3$ and $F$ is globally generated,  

(ii).- $h^0(F^{\vee})=0$.
\end{lem}

\noindent{\it{Proof.}} (i).- Consider $p_1,..., p_{g-5}$ general points on $C$. Note that  the line bundle $A:=K_C(-p_1-\cdots-p_{g-5})$ is of degree $g+3$, base-point-free with $h^0(C,A)=5$. $M_1:=K_C\otimes A^{\vee}$ is a line bundle with only one section and $h^0(C,M_1(-q_i))=0, i=1, 2, 3$. The line bundle $M_2:=A(-D)$ is such that $|M_2|=g^1_{g}$ is base-point-free. Consider the extension map $\text{Ext}^1(M_2,M_1)\xrightarrow{\beta}\text{Hom }(H^0(M_2), H^1(M_1))$ which sends an extension 

$$((e): 0\to M_1\to E\to M_2\to 0)$$ 

\noindent to the coboundary map 

$$(\partial_e:H^0(M_2)\to H^1(M_1)=H^0(A)^{\vee}).$$ 

\noindent Any non-trivial extension in $\text{Ker}(\beta)$ corresponds to a bundle $F$ satisfying $h^0(F)=3$ and $\text{det}(F)=\wedge^2F=K_C(-D)$. We have the isomorphisms $\text{Ext}^1(M_2, M_1)\simeq H^1(M_1\otimes M_2^{\vee})\simeq H^0(A^{\otimes 2}(-D))^{\vee}$ and $\text{Hom}(H^0(M_2), H^1(M_1))\simeq H^0(M_2)^{\vee}\otimes H^0(A)^{\vee}$.

\vskip1mm

\noindent We prove that $\text{Ker}(\beta)$ is a vector space of dimension $g-5$: the map $\beta$ is dual to the multiplication map $m_D:H^0(M_2)\otimes H^0(A)\to H^0(M_2\otimes A)=H^0(A^{\otimes 2}(-D))$. By the base point free-pencil trick  applied to $M_2$, $\text{Ker}(m_D)=H^0(C,\mathcal O_C(D))$, then $\text{dim Ker}(m_D)=1$ and $m_D$ has cokernel of dimension $h^0(A^{\otimes 2}(-D))-9=g-5=\text{dim Ker}(\beta)$. 

\vskip1mm

\noindent Now consider an extension $((e): 0\to M_1\to F\to M_2\to 0)$ in $\text{Ker}(\beta)$ and suppose that $F$ is not globally generated, then the three sections of $F$ generates a subsheaf $F_1$ of $F$ fitting into an exact sequence $0\to M_1(-B)\to F_1\to M_2\to 0$, for some divisor $B$ over $C$. For any $y\in C-\{p_1,...,p_{g-5}\}$, $h^0(M_1(-y))=0$ and $h^0(M_2(-y))=1$, then $h^0(F(-y))\leq 1$. By Riemann-Roch we have that $h^0(F)-h^0(F^{\vee}\otimes K_C)=-3$, that is, $h^0(K_C\otimes F^{\vee})=6$. This implies that $h^0(F(-y))=-5+h^0(K_C(y)\otimes F^{\vee})\geq 1$, so $h^0(F(-y))=1$ and $F$ is globally generated at $C-\{p_1,...,p_{g-5}\}$. This implies that $B\subset\text{supp}(\{p_1,...,p_{g-5}\})$. Thus, if $F$ comes from an element $e\in\text{Ker}({\beta})$ that fails to be generated by global sections, then there is a point $x$ among the points $\{p_1,...,p_{g-5}\}$ and a subsheaf $F_2$ given by the extension $0\to M_1(-x)\to F_2\to M_2\to 0$ so that the extension $(e)$ is induced from this latter extension and it is surjective on global sections. Since $\text{Hom}(H^0(M_2), H^1(M_1(-x)))\simeq\text{Hom}(H^0(M_2), H^0(A(x))$, such extensions are parametrized by elements in 

$$\text{Ker }[\text{Ext}^1(M_2,M_1(-x))\to\text{Hom}(H^0(M_2), H^0(A(x))],$$ 

\noindent where $\text{Ext}^1(M_2,M_1(-x))\simeq H^0(A^{\otimes 2}(x-D)))^{\vee}$. Note that since $h^0(A)=5$ and $x\in\text{supp }(\{p_j\}_{j=1}^{g-5})$ then $h^0(A(x))=6$. By the base point free pencil trick applied to $M_2$, note that $H^0(C,\mathcal O_C(x+D))$ is the kernel of the map $H^0(M_2)\otimes H^0(A(x))\to H^0(A^{\otimes 2}(x-D))$ and $h^0(C;\mathcal O_C(x+D))=1$, then the cokernel has dimension $h^0(A^{\otimes 2}(x-D))-11=g-6$. This implies that the extensions in $\text{Ker}(\beta)$ which fail to be generated by global sections have codimension at least 1 in 
$\text{Ker}(\beta)$, so for a general extension in $\text{Ker}(\beta)$, the corresponding vector bundle $F$ satisfies that $F$ is globally generated, $h^0(F)=3$ and $\wedge^2F=K_C(-D)$. 

\vskip1mm

\noindent (ii).- We have the identification $F\simeq F^{\vee}\otimes K_C(-D)$. Since $F$ is globally generated, one has: 

\begin{align}
0\to (K_C(D))^{\vee}\to H^0(C,F)\otimes\mathcal O_C\to F\to 0.
\end{align} 

\noindent Take $V=H^0(C,F)^{\vee}$ and dualize (1) to get 

\begin{align}0\to F^{\vee}\to V\otimes\mathcal O_C\to K_C(-D)\to 0,
\end{align} 

\noindent thus $h^0(F^{\vee}\otimes K_C(-D))=h^0(F)=3$. If $0\neq s\in H^0(F^{\vee})$, then we have an injection $s:\mathcal O_C\hookrightarrow F^{\vee}$, so $K_C(-D)\hookrightarrow F^{\vee}\otimes K_C(-D)=F$ and $3=h^0(F)\geq h^0(K_C(-D))=g-3$ which is a contradiction for $g\geq 7$, then $h^0(F^{\vee})=0$. $\square$
 
 \vskip1mm
 
\begin{lem} Let $\Bbb G:=G(3,H^0(C,K_C(-D)))$ be the Grassmannian of $3-$planes in $H^0(C,K_C(-D))$. For $V\in\Bbb G$ general, the map 

$$\mu_V:V\otimes H^0(C,K_C(-D))\to H^0((K_C(-D))^{\otimes 2})$$ 

\noindent as in (+) has kernel of dimension 3.\end{lem}

\vskip1mm

\noindent{\it{Proof.}}  Note that $\wedge^2V\subset\text{Ker}(\mu_V)$, then $\text{dim Ker} (\mu_V)\geq 3, \forall\hskip1mm V\in\Bbb G$. Let $\Sigma_3:=\{V\in\Bbb G:\text{dim Ker}(\mu_V)=3\}$. We are going to show that $\Sigma_3\neq\emptyset$ and that $\text{dim}(\Sigma_3)=\text{dim}(\Bbb G).$  

\vskip1mm

\noindent  Consider the vector bundle $F$ constructed in Lemma 3.1. Tensoring (2) by $K_C(-D)$ we have 

\begin{align}
0\to F^{\vee}\otimes K_C(-D)\to H^0(F)^{\vee}\otimes K_C(-D)\to (K_C(-D))^{\otimes 2}\to 0.
\end{align}

\noindent We know that $h^0(F^{\vee})=0$, also a non-zero section $\tau\in H^0(\mathcal O_C)$ induces an isomorphism $H^0(F)^{\vee}\otimes H^0(\mathcal O_C)\simeq H^0(F)^{\vee}$, hence taking cohomology in (2) we have an injective map $\iota: H^0(F)^{\vee}\hookrightarrow H^0(K_C(-D))$. Set $V:=\iota(H^0(F)^{\vee})$, then $V\subset H^0(K_C(-D))$ has dimension 3. Since $F^{\vee}\otimes\text{det}(F)\simeq F$, then $H^0(F^{\vee}\otimes K_C(-D))\simeq H^0(F)$, so we take cohomology  in (3) to obtain
\begin{align}
0\to H^0(C,F)\to H^0(F)^{\vee}\otimes H^0(K_C(-D))\to H^0((K_C(-D))^{\otimes 2})\to\cdots.
\end{align}
\noindent From (4) we have that  
\begin{align}
H^0(C,F)\simeq\text{Kernel}(H^0(F)^{\vee}\otimes H^0(K_C(-D))\to H^0((K_C(-D))^{\otimes 2})).
\end{align} 
\noindent Since $H^0(F)^{\vee}\otimes H^0(K_C(-D))\simeq V\otimes H^0(K_C(-D))$, from (5) we have $H^0(F)\simeq\text{Ker}(\mu_V)$, then $V\simeq H^0(C,F)$, so $\Sigma_3\neq\emptyset$. By upper semicontinuity of the function $\Bbb G\to\Bbb Z, V\to\text{dim Ker}(\mu_V)$, we have that for the general $V\in\Bbb G$, $\text{dim Ker}(\mu_V)=3$, then $\text{dim }(\Sigma_3)=\text{dim}(\Bbb G)$. $\square$

\vskip1mm

\noindent{\bf{Proof of the Theorem:}}

\vskip1mm

\noindent {\bf{Step 1.-}} Here we use the same strategy as in proofs of Theorems 5.8 and 5.17 in \cite{CF2}. Let $\Bbb P=\Bbb P((H^0(K_C(-D)^{\otimes 2})^{\vee})\simeq\Bbb P(\text{Ext}^1(K_C(-D), \mathcal O_C(D)))$. By Lemma 3.2 and the description of 
$\mathcal W_t$ as in ($\star$) we have that $\emptyset\neq\mathcal W_3\subset\Bbb P$ so it has expected dimension $3g-18$. Denote by $\pi_u$ the hyperplane in $\Bbb P$ defined by $\{u=0\}\subset H^0((K_C(-D))^{\otimes 2})$, where $u$ corresponds to the extension $(u): 0\to\mathcal O_C(D)\to E\to K_C(-D)\to 0.$

\vskip1mm

\noindent Let $\mathcal J_{\Bbb G}:=\{(W,\pi)\in\Bbb G\times\Bbb P:\text{Im}(\mu_W)\subset\pi\}$ be and let $\pi_1: \mathcal J_{\Bbb G}\to\Bbb G,\pi_2:\mathcal J_{\Bbb G}\to\Bbb P$ be the projections to the first and second factor respectively. From Lemma 3.2 we have that for $W\in\Bbb G$ general, $\text{dim(Im}(\mu_W))=3g-12$. The fiber of $\pi_1$ over a general element $V\in\Sigma_3$, is isomorphic to the linear system of hyperplanes in $\Bbb P$ passing through the linear subspace $\Bbb P(Im(\mu_V))$, then the general fiber of $\pi_1$ is irreducible and of dimension $(h^0((K_C(-D))^2)-1)-(3g-12)=2$. On the other hand, since 
$\Sigma_3$ is dense in $\Bbb G$ and $\Sigma_3\subseteq\pi_1(\mathcal J_{\Bbb G})$, then $\mathcal J_{\Bbb G}$ dominates 
$\Bbb G$ through $\pi_1$, thus there exists a unique component $\mathcal J_3$ dominating $\mathbb{G}$ through $\pi_1$ and $\text{dim}(\mathcal J_3)=2+\text{dim }(\Bbb G)=2+3(g-6)=3g-16$. By Serre duality, $\partial_u:H^0(C,K_C(-D))\to H^1(\mathcal O_C(D))$ is symmetric, that is, 
$\partial_u=\partial_u^{\vee}$, then $\text{Ker}(\partial_u)=\text{Ker}(\partial_u^{\vee})=(\text{Im}(\partial_u))^{\perp}$, in particular $\text{Ker}(\partial _u)$ is uniquely determined, so the general fiber of the map $\pi_2|_{\mathcal J_3}:\mathcal J_3\to\Bbb P$ is irreducible and zero-dimensional, then $\overline{\pi_2(\mathcal J_3)}\subset\Bbb P(\mathcal W_3)\subset\Bbb P $ is irreducible of dimension $3g-16$. Note that $\overline{\pi_2(\mathcal J_3)}$ gives rise to the existence of a component $\Delta_3\subset\mathcal W_3$ of dimension $3g-15$, which is therefore not good (in the sense of Definition 2.2) with $\Bbb P(\Delta_3)=\overline{\pi_2(\mathcal J_3)}$, such that for the general element $u\in\Delta_3$ we have that $\partial_u: H^0(C,K_C(-D))\to H^1(C,\mathcal O_C(D))$ has cokernel of dimension 3.

\vskip1mm

\noindent Consider the image $X$ of the map $C\hookrightarrow\Bbb P$ defined by the (very ample) linear system $|(K_C(-D))^{\otimes 2}|$. Note that for $\epsilon\in\{0,1\}$ and for $\sigma=g-6-\epsilon>0$ we have that $3g-16=\text{dim }(\Bbb P(\Delta_3))>\text{dim }\text{Sec}_{\frac{1}{2}(3g-16-\epsilon)}(X)$, then by Proposition 2.1 we have that for general $u\in\Delta_3$, the associated bundle $E_u$ is stable. This shows in particular that $B^4(C,K_C)\neq\emptyset$ and describes some points in it.

\vskip1mm

\noindent{\bf{Step 2.- Unobstructed sections.}} Let $u\in\Delta_3$ be a general extension and let $E=E_u$ be the associated vector bundle which fits into the exact sequence $(u): 0\to\mathcal O_C(D)\to E\to K_C(-D)\to 0$. Let $\Gamma$ be the section corresponding to the quotient $E\twoheadrightarrow K_C(-D)$. Let $S=\Bbb P(E)$ and $p:S\to C$ be the structure map. Let $c$ be the class of $\mathcal O_S(1)$ in $\text{Pic}(S)$ (or $\text{Num}(S)\simeq H^2(S,\Bbb Z))$. The tangent bundle $T_S$ fits in an exact sequence 
\begin{align}
0\to T_{S/C}\to T_S\to p^*T_C\to 0,
\end{align}  
\noindent where $T_{S/C}:=\text{Ker}(T_S\to p^*(T_C))$ is the relative tangent sheaf (the sheaf of tangent vectors along the fibers of $p$). The sheaf $T_{S/C}$ is dual to the relative canonical sheaf $\omega_{S/C}$ and $\omega_{S/C}=\mathcal O_S(-2c)\otimes p^*(\text{det}(E))=\mathcal O_S(-2c)\otimes p^*(K_C)$, then $T_{S/C}=\mathcal O_S(2c)\otimes p^*(T_C)$. On the other hand we have that $\mathcal N_{\Gamma/S}\simeq K_C(-2D)$, then $h^0(\Gamma, \mathcal N_{\Gamma/S})=g-6$ and $h^1(\Gamma, \mathcal N_{\Gamma/S})=1$. By (\cite{S}, p. 177, eq. (3.56)) We have 

\begin{align}
0\to T_S(-\Gamma)\to T_S<\Gamma>\to T_{\Gamma}\to 0.
\end{align}
\noindent Tensoring the exact sequence (6) by $\mathcal O_S(-\Gamma)$ we have
\begin{align}
0\rightarrow T_{S/C}(-\Gamma)\rightarrow T_S(-\Gamma)\rightarrow p^*(T_C)\otimes \mathcal O_S(-\Gamma)\rightarrow 0.
\end{align}

\noindent Recall that $\Gamma\sim\mathcal O_S(1) -f_D$ (see Section 2.1). Since $K_S\equiv -2c+(4g-4)f$, then $K_S(\Gamma)\equiv -c+a_0f$ for some $a_0\in\Bbb Z$. From the isomorphism $T_{S/C}\simeq \mathcal{O}_S(2c)\otimes p^*(T_C)$ we have that $K_S(\Gamma)\otimes(T_{S/C})^{\vee}\equiv -3c+af$ for some integer $a$, in particular $-3c+af$ is non-effective, then $h^2(T_{S/C}(-\Gamma))=h^0(K_S(\Gamma)\otimes (T_{S/C})^{\vee})=0$. Similarly we have that $K_S(\Gamma)-p^*(T_C)\equiv-c+bf$ for some integer $b$; then it is non-effective, so 
 $h^2(p^*(T_C)\otimes \mathcal{O}_S(-\Gamma))=h^0(K_S(\Gamma)\otimes(p^*(T_C))^{\vee})=0$. From (8) we have that $h^2(T_S(-\Gamma))=0$, and from (7) we deduce that $h^2(T_S<\Gamma>)=0$, then $(\Gamma,S)$ is unobstructed (see Section 2.3), that is, the first-order infinitesimal deformations of the closed embedding $\Gamma\hookrightarrow S$ are unobstructed with $S$ not fixed, in particular $\Gamma$ is unobstructed in $S$ fixed and $\Gamma$ varies in a $(g-6)$-dimensional family.

\vskip1mm

\noindent {\bf{Step 3.- Specially isolated sections.}} Let $$(u): 0\to\mathcal O_C(D)\to E=E_u\to K_C(-D)\to 0$$ 
\noindent be a general extension in $\Bbb P(\Delta_3)$. Consider the quotient $E_u\to\!\!\to K_C(-D)$ and the corresponding section $\Gamma_u=\Gamma$. To show that $\Gamma$ is specially isolated, we need to show that there are only finitely many sections corresponding to special quotients like these, that is, we need to show that the family of such corresponding quotients (which vary in a $(g-6)$-dimensional family) does not intersect in positive dimension the $3$-dimensional family $\mathcal F=\{L_t\}_{t\in T\subset\text{Sym}^3(C)}$ of special line bundles. By construction, we can assume that the family $\mathcal F$ is (isomorphic to) $W^{g-4}_{2g-5}(C)\simeq W_3^0(C)$.

\vskip1mm

\noindent Associated to the extension $(u)$ we have an exact cohomology sequence  

\vskip1mm

\noindent $0\to H^0(C,\mathcal O_C(D))\to H^0(C,E)\xrightarrow{\alpha}H^0(C,K_C(-D))\xrightarrow{\partial_u}H^1(C,\mathcal O_C(D))\to\cdots$

\vskip1mm

\noindent where by Step 1, $\text{Im}(\alpha)=\text{Ker}(\partial_u)$ is a general element in $\Bbb G$. Suppose that $h^0(C,E(-D))>1$, then there exists a non-zero section $s\in H^0(C,E), s\notin H^0(C,\mathcal O_C(D))$ such that $s$ vanishes at the points of $D$, then we have that $\alpha(s)\in\text{Ker}(\partial _u)\subset H^0(C,K_C(-D))$ is a non-zero section that vanishes at the points of $D$, that is, $\alpha(s)\in\text{Ker}(\partial_u)\cap H^0(C,K_C(-2D))$, thus $\text{Ker}(\partial_u)$ is an element of the proper and closed sublocus $\{W\in\Bbb G :W\cap H^0(C,K_C(-2D))\neq\emptyset\}\subsetneq\Bbb G$, since $\text{Ker}(\partial_u)$ is general in $\Bbb G$ this contradiction shows that $h^0(C,E(-D))=1$, which implies by Section 2.1 that $\text{dim}|\mathcal O_S(\Gamma)|=0$, that is, $\Gamma$ is linearly isolated (see Definition 2.1). 

\vskip1mm

\noindent Let $N=\mathcal O_C(D)$; by construction we have that $[N]$ is a general point in $W_3(E)$ (see Section 2.2). Since $h^0(E(-D))=1$, the fiber $\pi_3^{-1}([N])$ 
of the map $\pi_3:Q_3(E)\to W_3(E)$ at $[N]$ is only one point (see Section 2.2). Tensoring by $N^{\vee}=\mathcal O_C(-D)$ the exact sequence induced by $(u)$, we have in cohomology the exact sequence 

$$0\to H^0(\mathcal O_C)\to H^0(E(-D))\to H^0(K_C(-2D))\xrightarrow{\delta}H^1(\mathcal O_C)\to\cdots$$

\noindent where the map $\delta$ can be identified with the differential of $\pi_3:Q_3(E)\to W_3(E)$, then 
$\delta$ is an isomorphism onto its image, that is, 
$$H^0(K_C(-2D))=T_{[N]}(Q_3(E))\simeq T_{[N]}(W_3(E))=\delta(H^0(K_C(-2D)))\subseteq H^1(C,\mathcal O_C).$$ 
\noindent Denote by $<, >$ the Serre duality pairing, and let $V_1:=H^0(K_C(-2D))$. With these identifications of tangent spaces, to prove that $\Gamma$ is specially isolated we need to prove that $T_{[N]}(W_3(E))\cap T_{[N]}(W^0_3(C))=(0)=\delta(V_1)\cap (\text{Im}(\mu_{N}))^{\perp}$, where for $j=1,2$, $\mu_{N^j}:H^0(\mathcal O_C(jD))\otimes H^0(K_C(-jD))\to H^0(C,K_C)$ is the Petri map. Since $\text{Im}(\mu_N)=H^0(K_C(-D))\supset\text{Im}(\mu_{N^2})=V_1$, then $(\text{Im}(\mu_N))^{\perp}\subset V_1^{\perp}$. Let $\delta(\omega)\in\delta(V_1)\cap(\text{Im}(\mu_{N}))^{\perp}\subset\delta(V_1)\cap V_1^{\perp}$. Then $<\delta(\omega),v>=0$ for all $v\in V_1$, that is, $\delta(\omega)\in\text{Ker}(V_1\xrightarrow{\beta}V_1^{\vee})$, where $\beta: x\to\beta_x, \beta_x(v)=<x,v>$. By duality, $\beta$ is an isomorphism, then $\delta(\omega)=0$; since $\delta$ is injective, $\omega=0$. This prove that $\Gamma$ is specially isolated.

\vskip1mm

\noindent {\bf{Step 4. The map $\Bbb P(\Delta_3)\to B^4(2,K_C)$ is generically injective.}} By what proved in Step 1 we can consider the map $\pi: \Bbb P(\Delta_3)\to B^4(2,K_C)$
\begin{align}
((u): 0\to\mathcal O_C(D)\to E_u\to K_C(-D)\to 0)\longrightarrow [E_u].\end{align}
\noindent For a general element $[E_u]\in\pi(\Bbb P(\Delta_3))\subset B^4(2,K_C)$, we have that $\pi^{-1}([E_u])$ corresponds to the extensions $u'\in\Bbb  P(\Delta_3)$ that induce exact sequences $0\to\mathcal O_C(D)\to E_{u'}\to K_C(-D)\to 0$ and a diagram
\begin{align}
\xymatrix{
0\ar[r] & \mathcal O_C(D)\ar[r]^{\iota_1} & E_u\ar[r]\ar[d]_{\phi} & K_C(-D)\ar[r] & 0\\
0\ar[r] & \mathcal O_C(D)\ar[r]^{\iota_2} & E_{u'} \ar[r] & K_C(-D)\ar[r]  & 0
 }\end{align}
\noindent with an isomorphism $\phi:E_u\to E_{u'}$ of stable bundles. The maps $\phi\circ\iota_1$ and $\iota_2$ determine two non-zero sections $\sigma_1\neq \sigma_2$ in $H^0(E_u(-D))$. For Step 3, we have that for the general extension $u\in\Delta_3$, the corresponding section $\Gamma $ of the quotient line bundle $E_u\twoheadrightarrow K_C(-D)$ is linearly isolated, then $h^0(C,E_u(-D))=1$, so there exists some scalar $\lambda\neq 0$  such that $\phi\circ\iota_=\lambda\iota_2$. Since $E_u$ is stable, by Lemma 4.5 in (\cite{CF2}) we have that $E_u=E_{u'}$, i.e $u,u'$ are proportional in 
$\Delta_3$, then $\pi$ is generically injective. In particular for $u,u'\in\Delta_3$ general points, the vector bundles $E_u, E_{u'}$ cannot be isomorphic.

\vskip1mm

\noindent {\bf{Step 5. A global moduli map.}} Given a special and effective line bundle $L\in\text{Pic }(C)$ which is assumed to be a quotient line bundle of $E$, the condition for $E$ to have canonical determinant forces the kernel $N$  of the quotient map $E\twoheadrightarrow L$ to be isomorphic to $K_C\otimes L^{\vee}$. For $D$ a general effective divisor of degree three, $L=K_C(-D)$ depends on $\rho(L):=\rho(g,g-4,2g-5)=3$ parameters. From Steps 1-3 and the construction in ([5], Section 6) we obtain an irreducible component $\Bbb P(\widetilde{\Delta}_3)\subset\Bbb P(\mathcal E)$ of dimension $3g-16+\rho(L)=3g-13$, where a point in $\Bbb P(\widetilde{\Delta_3})$ corresponds to the datum of a pair $(y,u)$ with $y=(\mathcal O_C(D),K_C(-D))$, $D$ a general and effective divisor of degree 3 and $u\in\Bbb P(\Delta_3)$ an extension as in Step 1. From Steps 3 and 4 we have that the global moduli map ${\pi}_{2g-2,2g-5}|_{\Bbb P(\widetilde{\Delta_3})}:\Bbb P(\widetilde{\Delta_3})\to B^4(2,K_C)$ is generically injective and its image fills up an irreducible component $\mathcal B\subset B^4(2,K_C)$ of dimension 
$\rho_{K_C}(2,4,g)=3g-13$. $\square$

\end{document}